\documentstyle{article}
\begin{document}
\begin{large}

\centerline{{\bf A NOTE ON ANALYTICAL REPRESENTABILITY OF MAPPINGS}}

\centerline{{\bf INVERSE TO INTEGRAL OPERATORS}}

\centerline{{\bf M.I.OSTROVSKII}}

{\bf Abstract.} The condition onto pair ($F,G$) of function Banach spaces
under which there exists a integral operator $T:F\to G$ with analytic kernel
such that the inverse mapping $T^{-1}:$im$T\to F$ does not belong to 
arbitrary a priori given Borel (or Baire) class is found.
\medskip

We begin with recalling some definitions [4, \S 31.IX].  Let $X$
and $Y$ be metric spaces. By definition, the  family  of  {\it analytically
representable} mappings is the least family of mappings from $X$  to $Y$
containing all continuous mappings and being closed with respect  to
passage to pointwise limit.

This family is representable as an union
$\cup _{\alpha \in \Omega }\Phi _{\alpha }$, where $\Omega $ is the
set of all countable ordinals and $\Phi _{\alpha }$ are defined  in
the  following
way.

1. The class $\Phi _{0}$ is the set of all continuous mappings.

2. The class $\Phi _{\alpha } (\alpha >0)$ consists of all mappings
which are limits
of convergent sequences of mappings belonging to
$\cup _{\xi <\alpha }\Phi _{\xi }$.

A mapping $f:X\to Y$ is called of $\alpha $ {\it Borel class} if the set
$f^{-1}(F)$ is
a Borel set of multiplicative class $\alpha $ for every closed subset
$F\subset $Y.

It is known [1], [11], if $Y$ is a separable Banach  space,  then
the class $\Phi _{\alpha }$ coincides with the $\alpha $ Borel class for
finite $\alpha $ and  with
$\alpha +1$ Borel class for infinite $\alpha $. In addition [12],
if $Y$ is  a  Banach
space, then the class of all regularizable mappings  from $X$  to $Y$,
imortant in the theory of improperly posed problems coicides with $\Phi _{1}$.

The present paper is devoted to the following  problem.  Let $F$
and $G$ be Banach function spaces on [0,1] and $T:F\to G$ be  an  injective
integral operator with analytic kernel. What class  of  analytically
representable functions can the mapping $T^{-1}:TF\to F$ belong to?

Let us give some clarifications.

1. By analytic kernel we mean the mapping $K:[0,1]\times [0,1]\to {\bf C}$
with
the following property: For some open subsets $\Gamma $ and $\Delta $
of ${\bf C}$ such that
$[0,1]\subset \Gamma $ and $[0,1]\subset \Delta $, there exists
an analytic continuation of $K$  to
$\Gamma \times \Delta $.

2. If the spaces are real then we  consider  mapping $K$  taking
real values on $[0,1]\times [0,1]$.

The main result of the present paper states that for  wide
class of pairs $(F,G)$ and every countable ordinal $\alpha $ there  exists  an
injective integral operator $T:F\to G$ with an analytic kernel such  that
$T^{-1}:TF\to F$ does not belong to $\Phi _{\alpha }$.

This result is a generalization of the  Menikhes'  result  [6].
The last states that there exists an integral operator  from $C(0,1)$
to $L_{2}(0,1)$ with infinite differentiable kernel and  nonregularizable
inverse. A.N.Plichko [9]  generalized  this  result  onto  wide
classes of function spaces. The relations of our  results  with  the
results of [9] will be discussed in remark 3 below.

We use standard Banach space terminology and notation,  as  may
be found in [5].

For a subset $A$ of a Banach space $X$ by cl$A$, lin$A$ and $A^{\perp }$
we shall
denote the closure of $A$ in the strong topology, the linear span of $A$
and $\{x^*\in X^*:(\forall x\in A)(x^*(x)=0)\}$ respectively.
For a subset $A$ of $X$ by $A^{\top}$
we shall denote $\{x\in X:(\forall x^*\in A)(x^*(x)=0)\}$.
By $w^*-\lim$  we  shall  denote
the limit in the weak$^*$ topology.

In the present paper we restrict ourselves to the case  when $F$
is a separable Banach function space on the  closed  interval  [0,1]
continuously and injectively embedded into $L_1(0,1)$. It is clear that
classical separable function spaces satisfy this condition. The case
of nonseparable spaces continuously and  injectively  embedded  into
$L_1(0,1)$ is discussed in remark 4 below.

Let $\Gamma $ be a bounded open subset of the complex plane  such  that
$\Gamma \supset [0,1]$. Let $f:\Gamma \to {\bf C}$ be an analytic function
continuous on the closure
of $\Gamma$ . (If we consider real spaces then we assume in addition that $f$
takes real  values  on  the  real  axis.)  Function $f$  generates  a
continuous functional on $L_1(0,1)$ by means of the formula
$$(f,x)=\int^1_0f(t)x(t)dt.$$

Since $F$ is continuously embedded into $L_1(0,1)$  then $f$  generates  a
continuous functional on $F$. Let us denote by $U$  the  subset  of $F^*$
consisting of all functionals  of such type. Let $M$=cl$U$. It is  clear
that $M$ is a closed subspace of $F^*$.

Every element of $F$ may be considered as a functional on $M$.  The
corresponding embedding of $F$ into $M^*$ we denote by $H$.

{\bf Remark 1.} The subspace $M$ is a total subspace of $F^*$. It  follows
from the following facts: 1) $F$ is injectively embedded into $L_1(0,1)$;
2) the set of functionals  generated  by  polynomials  is  total  in
$(L_1(0,1))^*$.

The main result of the present paper is the following

{\bf THEOREM 1.} Let cl$(H(F))$ be of infinite codimension  in $M^*$  and
let $G$ contains a linearly independent infinite sequence of functions
which have analytic continuations onto some open  region $\Delta $  of  the
complex plane such that $\Delta \supset [0,1]$.
Then for every countable ordinal $\alpha $
there exists a linear continuous injective integral  operator $T:F\to G$
with analytic kernel such that $T^{-1}:TF\to F$ does not
belong to $\Phi _{\alpha }$.

Let us introduce necessary  definitions.  Let $X$  be  a  Banach
space. {\it Weak$^*$ sequential closure} of a subset $V$ of $X^*$
is defined to be
the set of all limits of weak$^*$ convergent sequences  from  $V$.  Weak$^*$
sequential closure will be denoted by $V_{(1)}$.  For  ordinal $\alpha $
{\it weak$^*$
sequential closure of order} $\alpha $ of a subset $V$ of
$X^*$ is defined by  the
equality $V_{(\alpha )}=\cup _{\beta <\alpha }(V_{(\beta )})_{(1)}$.

{\bf THEOREM 2} (A.N.Plichko [10]). Let $X$ be a separable Banach space
and $Y$ be an arbitrary Banach space. Let $T:X\to Y$ be a continuous linear
injective operator. Let $T^{-1}:TX\to X$ belongs to $\Phi _{\alpha }$.
If $\alpha $ is finite  then
$(T^*Y^*)_{(\alpha )}=X^*$. If $\alpha $ is infinite then
$(T^*Y^*)_{(\alpha +2)}=X^*$. Conversely, if $\alpha $
is finite and $(T^*Y^*)_{(\alpha )}=X^*$ or if $\alpha $ is
infinite  and $(T^*Y^*)_{(\alpha +1)}=X^*$,
then the mapping $T^{-1}:TX\to X$ belongs to $\Phi _{\alpha }$.

The proof of Theorem 1 is based on Theorem 2 and the  following
proposition.

{\bf Proposition.} Let $F$ and $M$ satisfy the assumptions of Theorem  1.
Then for every countable ordinal $\beta $ there exists a subspace $K$  of $M$
such that $K_{(\beta )}\neq F^*$.

Proof. Using the result due to Davis and Johnson [2, p.~360] we
find a weak$^*$ null sequence $\{u_n\}^{\infty }_{n=1}$  in $M$  and
a  bounded  sequence
$\{v_{k}\}^{\infty }_{k=1}$ in $F^{**}$ such that for some partition
$\{I_{k}\}^{\infty }_{k=1}$ of  the  set  of
natural numbers onto the pairwise disjoint infinite subsets we shall
have
$$
v_k(u_n)=\cases{1,&if $n\in I_k$;\cr
0,&if $n\not\in I_k$.\cr}
$$

Using the arguments from the proof of Theorem III.1 in  [3]  we
select a weak$^*$ basic subsequence
$\{u_{n(i)}\}^{\infty }_{i=1}\subset \{u_n\}^{\infty }_{n=1}$  such  that  the
intersection $I_k\cap \{n(i)\}^{\infty }_{i=1}$ is infinite
for  every $k\in {\bf N}$.  Passing  if
necessary to a subsequence and renumbering
$\{u_{n(i)}\}^{\infty }_{i=1}$ we obtain  the
sequence $\{y_{j}\}^{\infty }_{j=0}$ such that
$$
v_k(y_{j})=\cases{1,&if $j$ can be represented in the form $j=n(n+1)/2+k$
with $k\le n$;\cr
0,&if it is not the case.\cr}
$$
Definition of the weak$^*$ basic sequence implies that  the  space
$Z=F/(\{y_{j}\}^{\infty }_{j=0})^{\top}$ has the basis
$\{z_{j}\}^{\infty }_{j=0}$ such that $y_{j}(z_k)=\delta _{j,k}$.  It  is
easy to see that the sequence $\{z_{j}\}^{\infty }_{j=0}$ is bounded away
from zero  and
that the set $\{\sum^k_{i=j}z_{i(i+1)/2+j}\}^{\infty  \infty }_{j=0,k=j}$
is bounded.

We may assume without loss of generality that $||z_{i}||\le 1$ for  every
$i\in {\bf N}$. Using the arguments from the proof of lemma 1 in [7] we find  a
subspace $N$ of cl(lin$\{y_{j}\}^{\infty }_{j=0})$
and a bounded sequence $\{h_n\}^{\infty }_{n=1}$  in $Z^{**}$
such  that the following conditions are satisfied:

a) If a weak$^*$ convergent sequence $\{x^*_{m}\}^{\infty }_{m=1}$
is contained in $N_{(\gamma )}$
for some $\gamma <\beta $ and $x^*=w^*-\lim_{m\to \infty }x^*_{m}$ then
$$
h_n(x^*)=\lim_{m\to \infty }h_n(x^*_{m})
$$
for every $n\in {\bf N}$.

b) There exists a collection
$\{x^*_{n,m}\}^{\infty  \infty }_{n=1,m=1}$ of vectors  of $N_{(\beta )}$
such that for every $k,n\in {\bf N}$ we have
$$
w^*-\lim_{m\to \infty }x^*_{n,m}=0;
$$
$$
(\forall m\in {\bf N})(h^*_k(x^*_{n,m})=\delta _{k,n}).
$$
Let $\{s^*_k\}^{\infty }_{k=1}$ be a total (over $F$)
sequence in $M$. Let $c_1=\sup _n||h_n||$.
Let $\nu _n>0\ (n\in {\bf N})$ be such that
$\sum^{\infty }_{n=1}\nu _n<1/(2c_1)$.

We shall identify $Z^*$ with its natural image in $F^*$. Without loss
of generality we may assume that $Z^*\neq F^*$.

Let us denote by $g_n (n\in {\bf N})$ norm-preserving extensions of $h_n$ onto
$F^*$. Let us introduce the operator $R:F^*\to F^*$ by the equality
$$
R(x^*)=x^*+\sum^{\infty }_{n=1}\nu _ng_n(x^*)s^*_n.
$$
It is easy  to  check  that $||R-I||\le 1/2$ (where $I$  is  the  identity
operator), hence $R$ is an isomorphism.

Let $K=R(N)$. In the same way as in [7] we prove that  for  every
$\gamma \le \beta $ we have
$$
K_{(\gamma )}=R(N_{(\gamma )}).
\eqno{(1)}$$
Since $R$ is  an  isomorphism  then  relations
$K_{(\beta )}=R(N_{(\beta )})$  and
$N_{(\beta )}\subset Z^*\neq F^*$ imply that
$K_{(\beta )}\neq F^*$. At the same time  (1)  and b)  imply
that $R(x^*_{n,m})=x^*_{n,m}+\nu _ns^*_n\in K_{(\beta )}$ and therefore
$s^*_n\in K_{(\beta +1)}$. By the totality
of $\{s^*_n\}$ it follows that the subspace $K\subset F^*$ is total.
The proof of the proposition is complete.

Proof of Theorem 1. Using the proposition for $\beta =\alpha +2$ we  find  a
subspace $K$ of $M$ such that $K_{(\alpha +2)}\neq F^*$.
Using well-known arguments (see
e.g. [5, p.~43,44]) we find a fundamental minimal sequence
$\{f_i\}^{\infty }_{i=1}$
in the space $F$ such that its biorthogonal  functionals
$\{f^*_i\}^{\infty }_{i=1}$  are
total and are contained in $K$. It is easy to see (see e.g. [8])  that
there exists an isomorphism $S:F\to F$ such that functionals
$S^*f^*_i (i\in {\bf N})$
are contained in $U$. Let us introduce the notation
$g^*_i=S^*f^*_i (i\in {\bf N})$.  It
is clear that
$$
(\hbox{cl(lin}(\{g^*_i\}^{\infty }_{i=1})))_{(\alpha +2)}\neq F^*.
\eqno{(2)}$$

Reducing if necessary the domain $\Delta $ we may assume that
$\Delta $  is  bounded
and that $G$ contains a linearly  independent  sequence  of  functions
which have analytic continuations onto $\Delta $ which are continuous on the
closure of $\Delta $. Using standard biorthogonalization procedure (see e.g.
[5, p.~43,44]) we find a minimal sequence $\{q_i\}^{\infty }_{i=1}$
in $G$  such  that
there exist analytic continuations $\{\tilde{q}_i\}^{\infty }_{i=1}$
of $\{q_i\}^{\infty }_{i=1}$
to $\Delta $ which are
continuous on the closure of $\Delta $.

Let us introduce numbers
$$
a_i=\max \{||q_i||_G,\sup _{t\in \Delta }| \tilde{q}_i(t)| \}, (i\in {\bf N}).
$$
Let $r_i (i\in {\bf N})$ be functions analytic in $\Gamma $
and continuous  on  the
closure of $\Gamma $ such that functionals generated by them on $F$  coincide
with $g^*_i (i\in {\bf N})$. Let us introduce the numbers
$$
b_i=\max \{||g^*_i||_F*,\sup _{t\in \Gamma }| r_i(t)| \}, (i\in {\bf N}).
$$
Let us introduce the operator $T:F\to G$ by the equality
$$
T(f)=\sum^{\infty }_{i=1}g^*_i(f)q_i/(2^ia_ib_i).
\eqno{(3)}$$
This operator may be considered as an  integral  operator  with
the kernel
$$
K(t_1,t_2)=\sum^{\infty }_{i=1}r_i(t_1)\tilde{q}_i(t_2)/(2^ia_ib_i).
$$
By the definitions of $a_i$ and $b_i$ it  follows  that  this  series
converges to a function analytic in $\Gamma \times \Delta $.

Operator $T$ is injective since the sequence $\{g^*_i\}^{\infty }_{i=1}$
is total and
the  sequence $\{q_i\}^{\infty }_{i=1}$  is  minimal.  It  is  easy
to   see   that
$T^*G^*\subset $cl(lin$(\{g^*_i\}^{\infty }_{i=1}))$. By (2) we obtain
$(T^*G^*)_{(\alpha +2)}\neq F^*$. By Theorem 2
it follows that $T^{-1}$ does not belong to $\Phi _{\alpha }$.
The theorem is proved.

{\bf Remark 2.} It is easy to see  that  if $F$  and $G$  are  infinite
dimensional function spaces  and $T:F\to G$  is  an  injective  integral
operator with an analytic kernel then $G$ satisfies the  condition  of
Theorem 1.

{\bf Remark 3.} If the subspace $M\subset F^*$ is  norming  then  the  infinite
codimension  of cl$(H(F))$  in $M^*$  is  equivalent  to  the  infinite
codimension of lin$(M^{\perp }\cup F)$ in $F^{**}$. Therefore
Theorems 1 and 2  of  [9]
may  be  obtained  using  the  arguments  of  the   present   paper.
Furthermore, it follows that many of the spaces  considered  in  [9]
satisfy the  condition  of  Theorem  1  of  the  present  paper.  In
particular it is the case for classical nonreflexive spaces.

{\bf Remark 4.} Let $F$ be a  nonseparable  function  space,  which  is
continuously injectively embedded into $L_1(0,1)$, and $G$ is  the  space
satisfying the condition of  Theorem  1.  Let $\{g^*_i\}^{\infty }_{i=1}$
be  a  total
sequence in $F^*$ such that functionals $g^*_i (i\in {\bf N})$  are
represented  by
functions analytic in the same open  subset  containing  [0,1].  Let
$\{q_i\}^{\infty }_{i=1}$ be the sequence introduced in the proof of Theorem 1.
Let us
introduce the operator $T:F\to G$ by (3). The operator $T$ is an  injective
integral operator with an analytic kernel.  It  is  clear  that  the
image of $T$ is separable.  Therefore  the  mapping $T^{-1}:TF\to F$  is  not
analytically representable since the image of  the  separable  space
under the action of analytically representable mapping is separable.

\centerline{REFERENCES}

1. Banach S. \"Uber analytisch darstellbare Operationen in  abstrakten
Raumen, Fund. Math. 17 (1931), 283--295.

2. Davis W.J., Johnson W.B. Basic sequences and norming subspaces in
non-quasi-reflexive Banach spaces, Israel  J.  Math.  14  (1973),
353--367.

3. Johnson W.B., Rosenthal H.P.  On $w^*$-basic  sequences  and  their
applications to the study  of  Banach  spaces,  Studia  Math.  43
(1972), 77--92.

 4. Kuratowski K. Topology, v. I, New York, Academic Press, 1966.

5. Lindenstrauss J., Tzafriri L.  Classical  Banach  spaces,  v.  I,
Berlin, Springer-Verlag, 1977.

6. Menikhes L.D. On the  regularizability  of  mappings  inverse  to
integral operators, Doklady AN SSSR, 241 (1978), no.  2,  282--285
(Russian).

7. Ostrovskii M.I. Total  subspaces  with  long  chains  of  nowhere
norming weak$^*$ sequential closures, Note Mat., to appear.

8. Ostrovskii M.I. Regularizability of inverse linear  operators  in
Banach spaces with bases, Sibirsk. Mat. Zh.  33  (1992),  no.  3,
123--130 (Russian).

9. Plichko A.N. Non-norming subspaces and  integral  operators  with
non-regulari- zable inverses, Sibirsk. Mat. Zh. 29 (1988),  no.  4,
208--211 (Russian).

10. Plichko A.N.  Weak$^*$  sequential  closures  and $B$-measurability,
private communication.

11. Rolewicz S. On  the  inversion  of  non-linear  transformations,
Studia Math. 17 (1958), 79--83.

12. Vinokurov V.A. Regularizability and analytical representability,
Doklady AN SSSR 220 (1975), no. 2, 269--272 (Russian).

\end{large}
\end{document}